\documentclass[12pt]{amsart}
\usepackage{amssymb}
\usepackage{amsmath}
\usepackage{amsthm}

\makeatletter
\def\LaTeX{\leavevmode L\raise.42ex
    \hbox{\kern-.3em\size{\sf@size}{0pt}\selectfont A}\kern-.15em\TeX}
\makeatother

\newcommand{\BibTeX}{{\rm B\kern-.05em{\sc
          i\kern-.025emb}\kern-.08em\TeX}}

\makeatletter
\def\@currentlabel{2.1}\label{e:dispaa}
\def\@currentlabel{2.21}\label{e:dispau}
\def\@currentlabel{2.22}\label{e:dispav}
\def\@currentlabel{2.23}\label{e:dispaw}
\def\@currentlabel{2.24}\label{e:dispax}
\def\theequation{\thesection.\@arabic\c@equation}
\makeatother

\renewcommand{\theequation}{\arabic{section}.\arabic{equation}}

\newtheorem{thm}{Theorem}[section]

\newtheorem{cor}[thm]{Corollary}

\theoremstyle{definition}

\begin{document}
\title[quasi-linear elliptic system]
{Dirichlet problems of a quasi-linear elliptic system}
\author{Gongbo Li}
\address{Institute of Physics and Mathematical Sciences, WuHan, P.R. China}
\email{ligb@wipm.ac.cn }
\author{Li Ma}
\address{Department of Mathematical Sciences, Tsinghua University,
Beijing 100084, P.R. China} \email{lma@math.tsinghua.edu.cn}
\subjclass{Primary 35J50; Secondary 35J55} \keywords{elliptic
system, direct method, Dirichlet problem}
\date{}
\def\baselinestretch{1}

\begin{abstract}
We discuss the Dirichlet problem of the quasi-linear elliptic
system
\begin{eqnarray*}
-e^{-f(U)}div(e^{f(U)}\bigtriangledown
U)+&\frac{1}{2}f'(U)|\bigtriangledown U|^2&=0, \;\;\;\; \mbox{in
$\Omega$}, \;\;\;\; \\
& U|_{\partial\Omega}&=\phi.
\end{eqnarray*}
Here $\Omega$ a smooth bounded domain in $R^n$, $f: R^N\to R$ is a
smooth function, $U:\Omega\to R^N$ is the unknown vector-valued
function, $\phi:\overline\Omega\to R^N$ is a given vector-valued
$C^2$ function, $f'$ is the gradient of the function $f$ with
respect to the variable $U$.  Such problems arise in population
dynamics and Differential Geometry. The difficulty of studying
this problem is that this nonlinear elliptic system does not fit
the usual growth condition in M.Giaquinta's book [G] and the
natural working space $H^1\cap L^{\infty}(\Omega)$ for the
corresponding Euler-Lagrange functional does not fit the usual
minimization or variational argument. We use the direct method on
a convex subset of  $H^1\cap L^{\infty}(\Omega)$ to overcome these
difficulties. Under a suitable assumption on the function $f$, we
prove that there is at least one solution to this problem. We also
give application of our result to the Dirichlet problem of
harmonic maps into the standard sphere
\end{abstract}
\maketitle
\baselineskip 18pt
\section{Introduction}
\setcounter{equation}{0}

  In the interesting paper, Bensoussan-Boccardo-Murat [BBM] studied a
  class of quasi-linear elliptic equations which include the
  following equation as special case:
  $$
-\Delta u+u|Du|^2=h(x),\;\; \mbox{in $\Omega$}\;\;\;\;
  u|_{\partial \Omega}=0. \eqno(*)$$
where $\Omega$ a smooth bounded domain in $R^n$ with $h\in
L^{\infty}()\Omega$.  Observe that equation (*) has a variational
structure for left side of (*). In fact, let $\rho(u)=-u^2/2$.
Then the equation can be written as
$$-e^{-\rho(u)}div(e^{\rho(u)}\bigtriangledown
u)=h(x), \;\; \mbox{in $\Omega$}, \;\;\;\; u|_{\partial \Omega}=0.
$$
Let $P(u)=\int_0^ue^{\rho(s)} ds$. Clearly $P=P(u)$ has its
derivative $P'(u)=e^{\rho(s)}>0$, and  has its inverse function
$u=G(v)$. Define $v=P(u)$. Then the equation can be reduced into
$$
-\Delta v=e^{-G(v)^2/2}h:=h(x,v), \;\; \mbox{in $\Omega$},
\;\;\;\; v|_{\partial \Omega}=1.
$$
In this case, the nonlinear term is in right side and it is very
tamed. The variational structure is
$$
J(v)=\frac{1}{2}\int_{\Omega}|Dv|^2-\int_{\Omega} H(x,v)
$$
where $H(x,v)=\int^vh(x,s)ds$. Using the mountain pass method or
variational techniques one can easily prove the existence of a
classical positive solution of the problem (*). This kind of idea
is usually called the Cole-Hopf transformation in the study of
parabolic equations. One may also obtain a Liouville type theorem
for bounded smooth solutions of the problem:
$$-e^{-\rho(u)}div(e^{\rho(u)}\bigtriangledown
u)=0, \;\; \mbox{in $R^n$.}
$$
If we consider the evolution equation of the following form:
$$
u_t-\Delta u+u|Du|^2=h(x),\;\; \mbox{in
$\Omega$},\;\;\;\;u|_{t=0}=\varphi,\;\;\;\;
  u|_{\partial \Omega}=0. \eqno(**)$$
We can also use the trick above and reduced this equation into
$$
v_t-\Delta v=e^{-G(v)^2/2}h,\;\; \mbox{in
$\Omega$},\;\;\;\;v|_{t=0}=P(\varphi), \;\;\;\;
  v|_{\partial \Omega}=1. $$
Hence one can also get a global existence of solutions of the
problem (**). However, this trick is not so good for quasi-linear
elliptic systems or parabolic systems. We can also consider (*) in
another way. Let $m(u)=-u^2$. Then the elliptic equation in (*)
can be reduced into
$$-e^{-m(u)}div(e^{m(u)}\bigtriangledown
u)+\frac{1}{2}m'(u)|\bigtriangledown u|^2=h(x)
$$
and the left side is the variational derivative of the functional
$$
j(u)=\int_{\Omega} e^{m(u)}|Du|^{2}dx.
$$
This tells us that we can use this functional to study the
corresponding Dirichlet boundary value problem:
$$-e^{-m(u)}div(e^{m(u)}\bigtriangledown
u)+\frac{1}{2}m'(u)|\bigtriangledown u|^2=0,\;\; \mbox{in
$\Omega$}\;\;\;\;
  u|_{\partial \Omega}=\varphi.
$$
Although the functional $j(\cdot)$ is not so nice, one can show
that this problem is always solvable.

 In this paper, we mainly
discuss the Dirichlet problem of the quasi-linear elliptic system
$$-e^{-f(U)}div(e^{f(U)}\bigtriangledown
U)+\frac{1}{2}f'(U)|\bigtriangledown U|^2=0, \;\; \mbox{in
$\Omega$}, \;\;\;\; U|_{\partial \Omega}=\phi. \eqno(1.1) $$ Here
again $\Omega$ a smooth bounded domain in $R^n$, $f: R^N\to R$ is
a smooth function, $U:\Omega\to R^N$ is the unknown vector-valued
function, $\phi:\partial\Omega\to R^N$ is a given vector-valued
$C^{2, \delta}$ function, $f'$ is the gradient of the function $f$
with respect to the variable $U$.  Such problems arise in
population dynamics and Differential Geometry (see below). The
difficulty of this problem is that this nonlinear elliptic system
does not fit the usual growth condition in M.Giaquinta's book [G]
(see also [C] and [S]) and the natural working space $H^1\cap
L^{\infty}(\Omega)$ for the corresponding Euler-Lagrange
functional does not fit the variational argument. We point out
that the space $H^1\cap L^{\infty}(\Omega)$ or $W^{1,q}$ ($q>2$)
is the right for the corresponding variational functional to be
differentiable. Once we have this, we can derive the
Euler-Lagrange system for the functional. In our earlier research
[M], we noticed that there is one nice direction for us to
differentiate the variational functional. Then we use this
property and the Nash-De Giorgi- Moser-Morrey iteration to obtain
such a $L^{\infty}$ bound. This method is tricky, we used it  in
the study of certain elliptic systems (see [L]). However, it is
not so powerful. In fact, we can not use it to treat the system
(1.1).
 Our idea here is to find a
nice space obtaining a $L^{\infty}$ bound. To get such a
$L^{\infty}$ bound, we use the direct method on a convex subset of
$H^1\cap L^{\infty}(\Omega)$. In the following we will study the
existence of a weak solution of (1.1) in the class $H^1\cap
L^{\infty}(\Omega)$.

\begin{thm}
Suppose that $f(U)$ satisfies that $f'(U)=-Ug(U)$, where $g(U)$ is
a positive continuous function on $R^N$. Then there is at least
one weak bounded solution $U$ of (1.1). Furthermore, there is some
$q>2$ such that $U\in W^{1,q}(\Omega)$.
\end{thm}

Typical examples satisfying our assumption of Theorem 1 are (1)
$f(U)=-\alpha |U|^2$,(2) $f(U)=-2\log \frac{1+|U|^2}{2}$, and (3)
$f(U)=-\beta\log \frac{1+|U|^2}{2}$ where $\alpha, \beta>0$. Case
(1) corresponds to Gaussian measure $dm=e^{-|Y|^2}dY$. In this
case, $f'(Y)=-2\alpha Y$. Case (2) is for the standard metric
$ds^2=\frac{4}{(1+|Y|^2)^2} dY^2$on $S^{N-1}$ written in the
stereographic projection coordinates $(Y)$ from the north pole.
For this case, we have that $f'(Y)=-4Y/(1+|Y|^2)$. So the critical
map of $E(\cdot)$ is in fact a harmonic map from $\Omega$ to the
sphere. Then the result was obtained by S.Hildebrandt, H.Kaul, and
K.O.Widman [HKW], Schoen and Uhlenbeck in [SU] and Giaquinta and
Giusti in [GG]. When our $f$ satisfies some convexity condition,
one may consult the papers of S. Hildebrandt [H], D.G.Defigueiredo
[D], and Marcellini-Sbordone, [MS]. As a by-product, we can obtain
the following extension of a result due to V.Benci and J.M.Coron
in dimension two [BC].

\begin{cor}. Given a bounded domain $\Omega$ in $R^n$ with regular (Lipschitz)
boundary. Assume we have a $C^{2,\delta}$ map
$\phi:\partial\Omega\to S^{N}$. Auume $n\leq N$. Then these exist
at least two weakly harmonic maps from $\Omega \to S^{N}$ with
Dirichlet boundary value $\phi$.
\end{cor}
\begin{proof} Since $n\leq N$ and $\phi:\partial\Omega\to S^N$
is smooth, we have find a point $P$ in $S^N$ such that $\{P, -P\}$
are not in the range $\phi(\partial\Omega)$. Taking $P$ as the
north pole and using the stereographic coordinates at $P$, we can
reduce the harmonic map problem into our problem (1.1) with
$f(U)=-2\log \frac{1+|U|^2}{2}$. By Theorem 1.1, we have a weak
bounded solution $U$, which corresponds a weakly harmonic map from
$\Omega \to S^{N}$ with Dirichlet boundary value $\phi$. In the
same way, using $-P$ as south pole, we can obtain another weakly
harmonic map from $\Omega \to S^{N}$ with Dirichlet boundary value
$\phi$.
\end{proof}

Our method can be used to handle the following Dirichlet problem
of the quasi-linear elliptic system
$$-e^{-f(U)}\partial_i(e^{f(U)}A(x)_{ij}^{ab}\partial_j
U^a)+\frac{1}{2}[f'(U)]^b A(x)_{ij}^{ac}\partial_i U^a\partial_j
U^c=0, \;\; \mbox{in $\Omega$}, \;\;\;\; U|_{\partial
\Omega}=\phi. \eqno(1.1')
$$
Here we use the Einstein sum convention and we assume that
$(A^{ab}_{ij})$ is a uniformly positive matrix function in the
sense that its each component is non-negative and there are two
positive constants $\lambda$ and $\Lambda$ such that
$$
\lambda |\xi|^2\leq A(x)_{ij}^{ab}\xi^i_a\xi^j_b\leq \Lambda
|\xi|^2, \eqno(1.a)
$$
for any $\xi\in R^n\times R^N$. The energy integral for this
problem is
$$
   E(U)=\int_{\Omega} e^{f(U)}A^{ab}_{ij}\partial_iU^a\partial_jU^bdx
$$
Our result for this problem is

\begin{thm}
Suppose that $f(U)$ satisfies that $f'(U)=-Ug(U)$, where $g(U)$ is
a positive continuous function on $R^N$. Assume that the matrix
function $(A_{ij}^{ab}(x))$ satisfies $(1.a)$. Then there is at
least one weak bounded solution $U$ of (1.1'). Furthermore, there
is some $q>2$ such that $U\in W^{1,q}(\Omega)$.
\end{thm}

We also remark that we can extend our result to the Dirichlet
boundary problem on half space $R^n_+$
$$-e^{-f(U)}div(e^{f(U)}\bigtriangledown
U)+\frac{1}{2}f'(U)|\bigtriangledown U|^2=0, \;\; \mbox{in
$R^n_+$}, \;\;\;\; U|_{x_n=0}=\phi. \eqno(1.1'') $$ where $\phi$
is a smooth bounded function on $\overline{R^n_+}$.

 We have the following result:

\begin{thm}
Suppose that $f(U)$ satisfies that $f'(U)=-Ug(U)$, where $g(U)$ is
a positive continuous function on $R^N$. Assume that $\phi$ is a
smooth bounded function on $\overline{R^n_+}$ with
$\int_{R^n_+}|D\phi|^2dx<+\infty$. Then there is at least one weak
bounded solution $U$ of (1.1''). Furthermore, there is some $q>2$
such that $U\in W^{1,q}_{loc}(R^n_+)$.
\end{thm}

Clearly one can extend our results to the case when $f=f(U)$ is
replaced by a more general function $f=f(x,U)$. One may also
discuss the Dirichlet problem of a quasi-linear elliptic system
with a p-Laplacian type operator on the smooth bounded domain
$\Omega$ in $R^n$ or in the half space $R^n_+$:
$$-L_p U=0
\;\; \mbox{in $\Omega$}, \;\;\;\; U|_{\partial \Omega}=\phi.
\eqno(1.2)$$  Here, $p>1$,  and
$$L_pu:=e^{-f(U)}div(e^{f(U)}|DU|^{p-1}DU)+\frac{1}{2}f'(U)
|\bigtriangledown U|^2$$ with $f\in C^2(R^N)$.
 Since the formulations are complicated, we omit it in this paper.

\section{Proofs of Theorem 1.1 and 1.3}

 {\it We first prove Theorem 1.1}:

  We define the following integral
  $$
   E(U)=\int_{\Omega} e^{f(U)}|DU|^{2}dx
  $$
Then one can formally compute that
$$
<DE(U),V>=\int_{\Omega}e^{f(U)}(DU,
DV)+\frac{1}{2}\int_{\Omega}e^{f(U)}f'(U)|DU|^2dx
$$
Hence the Euler-Lagrange equation for $E$ is (1.1).

Since $\phi$ has a smooth extension over $\overline\Omega$, we can
find a constant vector $C$ such that $-C\leq \phi(x)\leq C$ for
each $x\in\overline\Omega$. Introduce the convex subset
$$
\mathcal{A}=\{U\in H^1(\Omega);-C\leq U(x)\leq C, \mbox{in
$\Omega$}, U(x)=\phi(x),\mbox{in $\partial\Omega$}\}
$$
Clearly $\phi\in \mathcal{A}$ and $\mathcal{A}$ is a weakly closed
convex subset of $H^1$. The functional $E$ is weakly lower
semi-continuous on $\mathcal{A}$ and $C^1$ on $H^1\cap
L^{\infty}(\Omega)$. Hence, we have a minimizer $U\in \mathcal{A}$
such that $E(U)=\inf_{v\in \mathcal{A}}E(v)$. We will prove that
this $U$ is a weak solution of (1.1). Then we prove Theorem 1.1.

Take any vector-valued function $\xi\in C^2_0(\Omega)$ and a small
positive constant $\epsilon$. Define
$$
V_{\epsilon}(x)=min\{C, max\{-C,U(x)+\epsilon\xi(x)\}\}.
$$
Here the maximum and minimum are taken for each component. Then we
have
$$
V_{\epsilon}\in\mathcal{A}.
$$
Define
$$
\Omega^{\epsilon}:=\{x\in \Omega; U(x)+\epsilon\xi(x)\geq C>U(x)\}
$$
and
$$
\Omega_{\epsilon}:=\{x\in \Omega; U(x)+\epsilon\xi(x)\leq
-C<U(x)\}
$$
Observe that for $\epsilon>0$ small we have $U>0$ on
$\Omega^{\epsilon}$ and $U<0$ on $\Omega_{\epsilon}$ . We also
have the measure $|\Omega^{\epsilon}|\to 0$ and the measure $
\Omega_{\epsilon}\to 0$ as $\epsilon\to 0+$.

Let
$$
\xi^{\epsilon}(x)=max\{0, -C+U(x)+\epsilon\xi(x)\}.
$$
and
$$
\xi_{\epsilon}(x)=-min\{0, C+U(x)+\epsilon\xi(x)\}.
$$
Clearly, $\xi^{\epsilon}(x)\geq 0$ and  $\xi_{\epsilon}(x)\geq 0$,
and both are in  $H^1_0\cap L^{\infty}(\Omega)$. Then we can
rewrite $V_{\epsilon}$ as
$$
V_{\epsilon}=U+\epsilon\xi-\xi^{\epsilon}+\xi_{\epsilon}.
$$

Since, for $t\in [0,1)$,
$$
E(U+t(V_{\epsilon}-U))\geq E(U),
$$
we have
$$
<DE(U),V_{\epsilon}-U>\geq 0.
$$
Note that
$$
<DE(U),V_{\epsilon}-U>=<DE(U),\epsilon\xi-\xi^{\epsilon}+\xi_{\epsilon}>.
$$
Hence we have that
$$
<DE(U),\xi>\geq
\frac{1}{\epsilon}[<DE(U),\xi^{\epsilon}>-<DE(U),\xi_{\epsilon}>]
$$
We will show that $<DE(U),\xi^{\epsilon}>\geq\circ (\epsilon)$ and
$<DE(U),\xi_{\epsilon}>\geq\circ (\epsilon)$. In fact, we have,
\begin{eqnarray*}
<DE(U), \xi^{\epsilon}>&=&\int_\Omega e^{f(U)}(DU,D\xi) dx\\
&+&\frac{1}{2} \int_\Omega (e^{f(U)}f'(U)\cdot \xi^{\epsilon} |DU|^2dx\\
&=& \int_{\Omega^{\epsilon}}
e^{f(U)}|DU|^2+\frac{1}{2}\int_{\Omega^{\epsilon}}
e^{f(U)}f'(U)\cdot
(U-C)|DU|^2\\
&+& \epsilon\int_{\Omega^{\epsilon}}
e^{f(U)}(DU,D\xi)+\frac{1}{2}\epsilon\int_{\Omega^{\epsilon}}
e^{f(U)}f'(U)\cdot
\xi|DU|^2\\
&=& \int_{\Omega^{\epsilon}}
e^{f(U)}|DU|^2g(U)(1-\frac{1}{2}U\cdot(U-C))\\
&+& \epsilon\int_{\Omega^{\epsilon}}
e^{f(U)}(DU,D\xi)+\frac{1}{2}\epsilon\int_{\Omega^{\epsilon}}
e^{f(U)}f'(U)\cdot \xi|DU|^2\\
 &\geq&
\epsilon\int_{\Omega^{\epsilon}}
e^{f(U)}(DU,D\xi)+\frac{1}{2}\epsilon\int_{\Omega^{\epsilon}}
e^{f(U)}f'(U)\cdot \xi|DU|^2=\circ(\epsilon).
\end{eqnarray*}
In the last inequality we used the fact that for $\epsilon>0$
small we have $U\cdot (C-U)>0$ on $\Omega^{\epsilon}$ (since $U>0$
and $C-U>0$ on $\Omega^{\epsilon}$), and in the last equality we
used $|\Omega^{\epsilon}|\to 0$. Similarly we can prove that
$<DE(U),\xi_{\epsilon}>\geq\circ (\epsilon)$.Therefore, we have
$$
<DE(U),\xi>\geq \circ(1)\to 0+
$$
Since our $\xi$ is arbitrary, we have
$$
<DE(U),\xi>= \circ(1)\to 0+
$$
and
$$
<DE(U),\xi>=0.
$$
Before finishing this section, we point out that our $U$ is a
spherical $Q-$ minima in the sense M.Giaquinta and Giusti [GG].
Hence we have $U\in W^{1,q}(\Omega)$ for some $q>2$. This proves
our Theorem 1.1.

The proof of Theorem 1.2 is a easy adaptation of the argument
above. So we omit the detail.

\section{Proof of Theorem 1.4}
\setcounter{equation}{0}

In this section we will give the proof of Theorem 1.3.

We choose a bounded domain exhaustion of $R^n_+$ in the way that
$R^n_+=U_{R>0} (R^n_+)\cap B_R(0)$. Here $B_R(0)$ is the ball of
radius $R$ with center at $0$. We write $B^+_R=(R^n_+)\cap
B_R(0)$.

For any bounded domain $\Omega$ in $R^n_+$, we let
$$
   E_{\Omega}(U)=\int_{\Omega} e^{f(U)}|DU|^{2}dx
  $$

Fix a large $R>0$.
 We define the following integral on $B^+_R$
  $$
   E_R(U)=\int_{B^+_R} e^{f(U)}|DU|^{2}dx
  $$
By Theorem 1.1, we have a bounded weak solution $u_R$ with uniform
bounds:
 $$ |u_R|_{L^{\infty}(B^+_R)}\leq |u_R|_{L^{\infty}(R_+^n)}
$$
and
$$ |Du_R|_{L^2(B^+_R)}\leq C_2 |D\phi|_{L^2(R_+^n))}
$$
where $C_2$ is a unform constant.

With these two bounds we can use extracting a diagonal subsequence
method to get a weakly convergence sequence $(U_k)$ and a weak
solution $W\in H^1_{loc}(R^n_+)\cap L^{\infty}(R^n_+)$ of
(1.1'')with the following property
$$
E_{\Omega}(W)\leq \underline{lim}_{k\to \infty} E_{\Omega}(U_k).
$$

This proves Theorem 1.3.

\end{document}